\documentclass[12pt]{article}
\usepackage{amsfonts,color}
\usepackage{amssymb}
\usepackage{amsmath}
\usepackage[dvips]{graphicx}

\newcounter{thm}[section]

\newtheorem{theorem}[thm]{Theorem}
\newtheorem{corollary}[thm]{Corollary}
\newtheorem{lemma}[thm]{Lemma}
\newtheorem{proposition}[thm]{Proposition}
\newenvironment{proof}[1][Proof]{\textbf{#1.} }{\ \rule{0.5em}{0.5em}}

\newcounter{assump}
\newtheorem{assum}[assump]{Assumption}

\begin{document}

\author{Olivier Lopez\thanks{%
Crest-Ensai and Irmar, rue Blaise Pascal, 35000 Bruz, France. E-mail : lopez@ensai.fr} \\
Crest-Ensai and Irmar}
\title{Single Index Regression Models with right censored responses}

\maketitle

\begin{abstract}
In this article, we propose some new generalizations of
M-estimation procedures for single-index regression models in
presence of randomly right-censored responses. We derive
consistency and asymptotic normality of our estimates. The results
are proved in order to be adapted to a wide range of techniques
used in a censored regression framework (e.g. synthetic data or
weighted least squares). As in the uncensored case, the estimator
of the single-index parameter is seen to have the same asymptotic
behavior as in a fully parametric scheme. We compare these new
estimators with those based on the average derivative technique of
Burke and Lu (2005) through a simulation study.
\end{abstract}


\bigskip

\section{Introduction}

In regression analysis, one investigates on the function $%
m(x)=E[Y\mid X=x]$, which is traditionally estimated from
independent copies $(Y_{i},X_{i})_{1\leq i\leq n}\in
\mathbf{R}^{1+d}$. The parametric approach consists of assuming
that the function $m$ belongs to some parametric family, that is
$m\left( x\right) =f_{0}\left( \theta _{0},x\right) $, where $%
f_{0} $ is a known function and $\theta_{0}$ an unknown finite
dimensional parameter. On the other hand, the nonparametric
approach requires fewer assumptions on the model, since it
consists of estimating $m$ without presuming the shape of the
function. However, this approach suffers from the so-called
''curse of dimensionality'', that is the difficulty to estimate
properly the function $m$ when the dimension $d$ is high (in
practice, $d\geq 3$). To avoid this important drawback of
nonparametric approaches, while allowing more flexibility than a
purely parametric model, one may use the semi-parametric
single-index model (SIM in the following) which states
\[
m\left( x\right) =E[Y\mid X^{\prime }\theta _{0}=x^{\prime }\theta
_{0}]=f\left( x^{\prime }\theta _{0};\theta _{0}\right) ,
\]
where $f$ is an unknown function and $\theta _{0}$ an unknown
finite dimensional parameter. If $\theta _{0}$ were known, the
problem would consist of a nonparametric one, but with the
covariates belonging nevertheless to a one-dimensional space.

In this framework, numerous semi-parametric approach have been
proposed for root-$n$ consistent estimation of $\theta _{0}$.
Typically, these approaches can be split into three mains
categories : $M$-estimation (Ichimura, 1993, Sherman, 1994b,
Delecroix et Hristache, 1999, Xia et Li, 1999, Xia, Tong, et Li,
1999, Delecroix, Hristache et Patilea, 2006), average derivative
based estimation (Powell, Stock et Stoker, 1989, H\"{a}rdle et Stoker,
1989, Hristache et al., 2001a, 2001b), and iterative methods
(Weisberg et Welsh, 1994, Chiou et M\"{u}ller, 1998, Bonneu et Gba,
1998, Xia et H\"{a}rdle, 2002).

If the responses of this regression model are randomly right-censored, these approaches clearly need to be adapted, for the
random variable $Y$ is not directly observed. The right
censoring model states that, instead of observing $Y$, one
observes i.i.d. replications of
\begin{eqnarray}
T &=&Y\wedge C, \nonumber \\
\delta &=&1_{ Y\leq C }, \label{censoring}
\end{eqnarray}
where $C$ is some ''censoring variable'', and $\mathbf{1}_A$
denotes the indicator function of the set $A$. In this setting,
semi-parametric Cox regression model (see e.g. Andersen et Gill,
1982) can be seen as a particular case of the SIM model, but
allows less flexibility. Moreover, it is still interesting to
extend mean-regression models to the censored framework. For this
reason, Buckley and James (1978) proposed an estimator of the
linear model under random censoring, and Lai and Ying (1991) and
Ritov (1990) proved its asymptotic normality. Koul, Susarla and
Van\ Ryzin (1981) initiated what we may call the ''synthetic
data'' approach, based on transformations of the data. See
Leurgans (1987), Zhou (1992b) and Lai \& al. (1995). Zhou (1992a)
also proposed a weighted least-square approach, applying weights
in the least square criterion in order to compensate the
censoring. These techniques were then used in the nonlinear
regression setting, that is when $f_{0}$ is known but nonlinear.
Stute (1999) established a connection between the weighted
least-square criterion and Kaplan-Meier integrals. Delecroix,
Lopez and Patilea (2006) extended the synthetic data approach.
Heuchenne and Van\ Keilegom (2005) modified the Buckley-James'
technique for polynomial regression purpose. When it comes to the
SIM model under random censoring, Burke and Lu (2005) recently
proposed an estimate using an extension of the average derivatives
technique of H\"{a}rdle and Stoker (1989) and the synthetic data
approach of Koul, Susarla, Van Ryzin (1981).

In this paper, we propose a semi-parametric $M$-estimator of the
SIM model under random censoring. We present a technique that is
adapted to both main classes of censored regression techniques
(synthetic data and weighted least squares), deriving root-$n$
consistency of our estimate of $\theta _{0}$, and then using it to
estimate $m\left( x\right) $. Another advantage of our technique
is that we do not require that the covariates $X$ have a density
with respect to Lebesgue's measure (only the linear combinations
$\theta'X$ need to be absolutely continuous), which is an important advantage comparatively
with the estimation procedure of Burke and Lu (2005).

The paper is organized as follows. In section \ref{secmodel} we
present the regression model and our methodology. In section
\ref{seccons}, we derive consistency of our semi-parametric
estimates in a general form, asymptotic normality is obtained in
section \ref{secas}. A simulation study is presented in
\ref{secsimul} to test the validity of our estimate with finite
samples. Section \ref{secappend} is devoted to technical
proofs.

\section{Model assumptions and methodology}

\label{secmodel}

In the following, we assume that we have the following regression model,
\[
Y=f\left( \theta _{0}^{\prime }X;\theta _{0}\right) +\varepsilon ,
\]
where $\theta_0$ is a vector of first component equal to $1,$ and
$E\left[ \varepsilon \mid X\right] =0$. The function $f$ is
defined in the following way, $
f\left( u;\theta \right) =E\left[ Y\mid X^{\prime }\theta =u\right] .
$
Considering the censoring model (\ref{censoring}), we will define
the following distribution function,
\begin{eqnarray*}
F\left( t\right) &=&\mathbb{P}\left( Y\leq t\right) , \\
G\left( t\right) &=&\mathbb{P}\left( C\leq t\right) , \\
H\left( t\right) &=&\mathbb{P}\left( T\leq t\right) , \\
F_{\left( X,Y\right) }\left( x,t\right) &=&\mathbb{P}\left( Y\leq
t,X\leq x\right) .
\end{eqnarray*}

In the following, we will assume that
\begin{eqnarray} \label{queue}
\inf\{t, F(t)=1\}&=&\inf\{t, H(t)=1 \}, \\
\mathbb{P}(Y=C) &=& 0. \label{pyc}
\end{eqnarray} Otherwise, if (\ref{queue}) does not hold,
since some part of the distribution of $Y$ remains unobserved,
consistent estimation requires making additional restrictive
assumptions on the law of the residuals. Note that, in this case, our estimators will still be
root-$n$ convergent, but not necessary to $\theta_0.$ Concerning (\ref{pyc}), we use this assumption to avoid dissymetry problems between
$C$ and $Y.$

As a property of conditional expectation, for any function $J(\cdot)\geq
0$, we have
\begin{eqnarray}
\theta _{0} &=&\arg \min_{\theta \in \Theta }E\left[ \left(
Y-f\left( \theta ^{\prime }X;\theta \right) \right)
^{2}J(X)\right] =\arg \min_{\theta \in \Theta
}M\left( \theta,f \right)  \label{M} \\
&=&\arg \min_{\theta \in \Theta }\int \left( y-f\left( \theta
^{\prime }x;\theta \right) \right) ^{2}J(x)dF_{(X,Y)}\left(
x,y\right). \nonumber
\end{eqnarray}
In equation (\ref{M}), of course we can not exactly know $\theta_0,$ since
two objects are missing in the definition of $M$, that is the distribution function
$F_{(X,Y)}$ and the regression function $f\left( \theta ^{\prime
}x;\theta \right) $. A natural way to proceed consists of
estimating these two functions, and then plugging in these
estimators into (\ref{M}).

\subsection{Estimating the distribution function}

We already mentioned there are two main approaches for studying
regression models in presence of censoring, the Weighted Least
Square approach (WLS in the following) and the Synthetic Data approach (SD in
the following).

\textbf{The WLS\ approach. } In the uncensored case, the distribution function $F_{(X,Y)}$ can be estimated using the empirical distribution. This tool is unavailable under random censoring, since it relies on the (unobserved) $(Y_i)_{1\leq i\leq n}.$
Under random censoring, Stute (1993)
proposed to use an estimator based on the Kaplan-Meier estimator of
$F$. Recall the definition of Kaplan-Meier estimator,
\[
\hat{F}\left( t\right) =1-\prod_{i:T_i \leq t}\left(
1-\frac{\sum_{j=1}^n \mathbf{1}_{\delta_j=1,T_j\leq
T_i}}{1-\hat{H}(T_i-)}\right)^{\delta_i},
\]
where $\hat{H}$ denotes the empirical distribution function of
$T$. $\hat{F}$ can be rewritten as 
\[
\hat{F}\left( y\right) =\sum W_{in}\mathbf{1}_{T_{i}\leq y},
\]
where $W_{in}$ is the jump at observation $i.$
It is particularly interesting to notice that the jump at observation $i$ is
connected to the Kaplan-Meier estimate of $G$ at the same value (see, for
example, Satten and Datta, 2000), that is
\begin{equation}
W_{in}=\frac{1}{n}\frac{\delta _{i}}{1-\hat{G}\left( T_{i}-\right) }.
\label{jump}
\end{equation}
Kaplan-Meier estimate is known to be a consistent estimate of $F$ under the two
following identifiability assumptions, that is
\begin{assum}
\label{a1} $Y$ and $C$ are independent.
\end{assum}
\begin{assum}
\label{a2} $\mathbb{P}\left( Y \leq C \mid
X,Y\right)=\mathbb{P}\left( Y \leq C \mid Y\right).$
\end{assum}

A major case for which Assumptions \ref{a1}-\ref{a2} hold is the case where
$C$ is independent from $\left( Y,X\right) $. However, Assumption
\ref{a2} is more general and covers a significant amount of
situations (see Stute, 1999).

\textbf{The SD approach. }
The SD approach consists of considering some alternative variable which has the same conditional
expectation as $Y$. 
For this, observe that, through elementary calculus, under Assumptions \ref{a1}-\ref{a2},
\begin{equation}
\forall \phi, \; E\left[ \frac{\delta \phi \left( X,T\right)
}{1-G\left( T-\right) }\mid X\right] =E\left[ \phi \left(
X,Y\right) \mid X\right]. \label{calcul}
\end{equation}
From (\ref{calcul}), we see that, if we
define, accordingly to Koul \& al. (1981),
\[
Y^{*}=\frac{\delta T}{1-G\left( T-\right) },
\]
we have $E\left[ Y^{*}\mid X\right] =E\left[ Y\mid X\right] $
under Assumption \ref{a1} and \ref{a2}. Hence, if $Y^{*}$ were
available, the same regressions techniques as in the uncensored
case could be applied to $Y^{*}$. Of course, $Y^{*}$ can not be
computed, since it depends on the unknown function $G$. But
$Y^{*}$ can be easily estimated (which is not the case for $Y$) by
replacing $G$ by its Kaplan-Meier estimate. For $i=1,...,n$ we
obtain
\[
\hat{Y}_{i}^{*}=\frac{\delta _{i}T_{i}}{1-\hat{G}\left( T_{i}-\right) }.
\]
See also Leurgans (1987), Lai \& al. (1995) for other kind of transformations.

Back to equation (\ref{M}), the SD approach will first consists of observing
that
\begin{eqnarray}
\theta _{0} &=&\arg \min_{\theta \in \Theta }E\left[ \left( Y^{*}-f\left(
\theta ^{\prime }x;\theta \right) \right) ^{2}J(X)\right]=M^*(\theta,f)  \label{bobo} \\
&=&\arg \min_{\theta \in \Theta }\int \left( y^{*}-f\left( \theta
^{\prime }x;\theta \right) \right)
^{2}J(x)dF_{(X,Y^{*})}^{*}\left( x,y^{*}\right) , \nonumber
\end{eqnarray}
where $F_{\left( X,Y^{*}\right) }^{*}\left( x,y^{*}\right)
=\mathbb{P}\left( X\leq x,Y^{*}\leq y^{*}\right) .$ 
%

Note that $M^*$ and $M$ are not the same functions. Indeed, $Y^*$
happens to have the same conditional expectation as $Y$ (hence $M$
and $M^*$ have the same minimizer $\theta_0$), but it has not the
same law.
\subsection{Estimating $f\left( \theta ^{\prime }x;\theta \right) $}

In the uncensored case, a common non-parametric way to estimate a
conditional expectation is to use kernel smoothing. In this case,
the Nadaraya-Watson estimate for $f\left( \theta ^{\prime
}x;\theta \right) $ is
\begin{eqnarray*}
\hat{f}\left( \theta ^{\prime }x;\theta \right) &=&\frac{\sum_{i=1}^{n}K%
\left( \frac{\theta ^{\prime }X_{i}-\theta ^{\prime }x}{h}\right) Y_{i}}{%
\sum_{i=1}^{n}K\left( \frac{\theta ^{\prime }X_{i}-\theta ^{\prime }x}{h}%
\right) } \\
&=&\frac{\int yK\left( \frac{\theta ^{\prime }u-\theta ^{\prime }x}{h}%
\right) d\hat{F}_{emp}\left( u,y\right) }{\int K\left( \frac{\theta ^{\prime
}u-\theta ^{\prime }x}{h}\right) d\hat{F}_{emp}\left( u,y\right) }.
\end{eqnarray*}
We are still facing the same problem of absence of the empirical
distribution function. However, WLS and SD approaches can be used
to extend the Nadaraya-Watson estimate to censored regression. In
the following, we will only use the SD approach of Koul \& al. to
estimate the conditional expectation, that is
\begin{equation}
\hat{f}\left( \theta ^{\prime }x;\theta \right) =\frac{\sum_{i=1}^{n}K\left(
\frac{\theta ^{\prime }X_{i}-\theta ^{\prime }x}{h}\right) \hat{Y}_{i}^{*}}{%
\sum_{i=1}^{n}K\left( \frac{\theta ^{\prime }X_{i}-\theta ^{\prime }x}{h}%
\right) }.  \label{ksvcond}
\end{equation}
While using this estimator, we do not have to deal with
Kaplan-Meier integrals at the denominator. In fact, the integral
at the denominator becomes an integral with respect to the
empirical distribution function of $X$. However, alternative
estimates (not necessarily kernel estimates) can still be used,
provided that they satisfy some further discussed conditions to
achieve asymptotic properties of $\hat{\theta}$. Therefore we
chose to present our results without presuming on the choice of
$\hat{f}\left( \theta ^{\prime }x;\theta \right) $, and then to
check in the Appendix section that the estimator defined in
(\ref{ksvcond}) satisfies the proper conditions.

Also observe that, using this kernel estimate, contrary to the
average derivative technique of Burke and Lu (2005), we do not
need to impose that $X$ has a density with respect to Lebesgue's
measure. We only need that the linear combinations $\theta'X$ do.

{\textbf{The choice of the trimming function $J$.}} The reason why
we introduced the function $J$ in (\ref{M}) appears in the
definition (\ref{ksvcond}). To ensure uniform consistency of this
estimate, we will need to bound the denominator away from zero.
For this, we will need to restrain the integration domain to a set
where $f_{\theta'X}(u)$ is bounded away from zero, $f_{\theta'X}$
denoting the density of $\theta'X.$ If we were to know $\theta_0,$
we could consider a set $B_0=\{ u : f_{\theta_0'X}(u)\geq c\}$ for
some constant $c>0$, and use the trimming
$J(\theta_0'X)=\mathbf{1}_{\theta_0'X \in B_0}.$ Of course,
this ideal trimming can not be computed, since it depends on the
unknown parameter $\theta_0$. Delecroix \& al. (2006) proposed a
way to approximate this trimming from the data. Given some
preliminary consistent estimator $\theta_n$ of $\theta_0,$ they use the following trimming,
$$J_n(\theta_n'X)=\mathbf{1}_{\hat{f}_{\theta_n'X}(\theta_n'X)\geq
c}.$$ In the following proofs, we will mostly focus on the
estimation using the uncomputable trimming $J(\theta_0'X)$, and we
will show in the appendix section that there is no asymptotic
difference in using $J_n(\theta_n'X)$ rather than $J(\theta_0'X)$.

\subsection{Estimation of the single-index parameter}

{\textbf{Preliminary estimate of $\theta_0$.}} For a preliminary
estimate, we assume, as in Delecroix \& al. (2006) that we know
some set $B$ such as $\inf_{x \in B, \theta \in
\Theta}\{f_{\theta'X}(\theta'x)\geq c>0 \},$ and we consider the
trimming function $\tilde{J}(x)=\mathbf{1}_{x \in B}$. To
compute our estimate $\theta_n,$ we then can use either of the WLS
or SD approach. For example, using the WLS approach, let
\begin{equation}
\theta_n=\arg \min_{\theta \in \Theta }\int \left( y-\hat{f}%
\left( \theta ^{\prime }x;\theta \right) \right)
^{2}\tilde{J}(x)d\hat{F}_{(X,Y)}\left( x,y\right)=\arg \min_{\theta \in \Theta }M_n^p(\theta,\hat{f}).
\label{preliminary}
\end{equation}

{\textbf{Estimation of $\theta_0$.}} In view of (\ref{M}) and
(\ref{bobo}), we will define our estimates of $\theta _{0}$
according to the two regression approaches discussed above,
\begin{eqnarray*}
\hat{\theta}_{WLS} &=&\arg \min_{\theta \in \Theta_n }\int \left[ y-\hat{f}%
\left( \theta ^{\prime }x;\theta \right) \right]
^{2}J_n(\theta_n'x)d\hat{F}_{(X,Y)}\left( x,y\right) \\
&=&\arg \min_{\theta \in \Theta_n }M_{n}^{WLS}\left( \theta,\hat{f} \right)
,
\\
\hat{\theta}_{SD} &=&\arg \min_{\theta \in \Theta_n }\int \left[ y^{*}-\hat{f}%
\left( \theta ^{\prime }x;\theta \right) \right]
^{2}J_n(\theta_n'x)d\hat{F}^{*}\left( x,y^{*}\right) \\
&=&\arg \min_{\theta \in \Theta_n }M_{n}^{SD}\left( \theta,\hat{f} \right) .
\end{eqnarray*}
In the definition above, for technical convenience, we restrained our optimization to
shrinking neighborhoods $\Theta_n$ of $\theta_0,$ chosen
accordingly to the preliminary estimation by $\theta_n.$ 

\subsection{Estimation of the regression function}
With at hand a root-n consistent estimate of $\theta_0$, it is
possible to estimate the regression function by using
$\hat{\theta}$ and some estimate $\hat{f}$. For example, using
$\hat{f}$ defined in (\ref{ksvcond}) will lead to
\begin{eqnarray*}
\hat{f}\left( \hat{\theta} ^{\prime }x;\theta \right)
=\frac{\sum_{i=1}^{n}K\left(
\frac{\theta ^{\prime }X_{i}-\hat{\theta} ^{\prime }x}{h}\right) \hat{Y}_{i}^{*}}{%
\sum_{i=1}^{n}K\left( \frac{\hat{\theta} ^{\prime }X_{i}-\hat{\theta} ^{\prime }x}{h}%
\right) }.
\end{eqnarray*}

\section{Consistent estimation of $\theta_0$}

\label{seccons}

In this section, we prove consistency of $\theta_n$ where $\theta_n$ is defined in (\ref{preliminary}). As a consequence, $\hat{\theta}$ is consistent since it is obtained from minimization other a shrinking neighborhood of $\theta_0.$ We will need two kinds of assumptions to
achieve consistency :  general assumptions on the regression model including identifiability assumptions for $\theta
_{0}$, and assumptions on $\hat{f}$.


\smallskip

\textbf{Identifiability assumptions for }$\theta _{0}$\textbf{ and
assumptions on the regression model.}

\begin{assum}
\label{a3} $EY^{2}<\infty $.
\end{assum}

\begin{assum}
\label{a4} If
$M(\theta_1,f)=M(\theta_0,f),$ then
$\theta_1=\theta_0$.
\end{assum}

\begin{assum}
\label{a6} $\Theta$ and $\mathcal{X}=Supp(X)$ are compact subsets
of $\mathbb{R}^d$ and $f$ is continuous with respect to $x$ and $\theta.$ Furthermore,
assume that $| f\left( \theta _{1}^{\prime }x;\theta
_{1}\right) -f\left( \theta _{2}^{\prime }x;\theta _{2}\right)|
\leq \left\| \theta _{1}-\theta _{2}\right\| ^{\gamma }\Phi \left(
X\right) $, for a bounded function $\Phi \left( X\right) $, and
for some $\gamma>0.$
\end{assum}

Assumption \ref{a3} is implicitly needed in order to define $M,$ while Assumption \ref{a4} ensures the identification of the parameter $\theta_0.$
On the other hand, Assumption \ref{a6} states that the class of functions $\mathcal{F=}\left\{
f\left( \theta ^{\prime }.;\theta \right) ,\theta \in \Theta
\right\} $ is sufficiently regular to allow it to satisfy an
uniform law of large numbers property. More precisely, Assumption \ref{a6} ensures that this class is Euclidean for a bounded envelope, according to Pakes and Pollard (1989). Observe that the condition that $\Phi $ is bounded can be
weakened, by replacing it by a moment assumption on $\Phi $.
However, this condition is quite natural in a context where we
will assume that the covariates are bounded random vectors, and
this will simplify our discussion. Moreover, it implies that $f$
is a bounded function of $\theta$ and $x.$


\textbf{Assumptions on }$\hat{f}.$

\begin{assum}
\label{a7} For all function $g$, define, for $c>0,$
$$\| g \|_{\infty}=\sup_{\theta \in \Theta,x} |g(\theta ' x;\theta )|\mathbf{1}_{f_{\theta}(\theta'x)>c/2}.$$
Assume that $\|\hat{f}-f \|_{\infty}=o_{P}(1).$
\end{assum}

See section \ref{secappend} for more details to see that the kernel estimator (\ref{ksvcond}) satisfies this assumption under some additional integrability assumptions on the variable $Y.$


\begin{theorem}
\label{consistency} Under Assumptions \ref{a3} to \ref{a7}, we
have
\[
\sup_{\theta \in \Theta }\left| M_{n}\left( \theta ,\hat{f}\right) -M_{\infty}\left( \theta,f \right)
\right| =o_{P}\left( 1\right).
\]
As an immediate corollary, in a probability sense, $\theta_n  \rightarrow  \theta _{0}$.
\end{theorem}

\begin{proof} 
%

\textbf{Step 1 : replacing $\hat{f}$ by $f.$} Observe that, since the integration domain is restricted to the set $B,$
\begin{eqnarray*}
\left|M_n(\theta,f)-M_n(\theta,\hat{f})\right| &\leq & \|\hat{f}-f\|_{\infty}\\&&\times[\|\hat{f}
 +f\|_{\infty}\int d\hat{F}_{(X,Y)}(x,y)+2\int |yd\hat{F}_{(X,Y)}(x,y)|].
\end{eqnarray*}
Now using Assumption \ref{a7}, deduce that $\sup_{\theta\in \Theta}|M_n(\theta,f)-M_n(\theta,\hat{f})|=o_P(1).$

\textbf{Step 2 : $M_n(\theta,f).$} Showing that $\sup_{\theta\in \Theta}|M_n(\theta,f)-M(\theta,f)|=o_P(1)$ can then be done in the same way as in a nonlinear regression model such as in Stute (1999). See the proof of Theorem 1.1 in Stute (1999).
\end{proof}

\section{ Asymptotic normality}

\label{secas}

As in the uncensored case, we will show that, asymptotically
speaking, our estimators behave as if the true family of functions $f$ were
known. Hence studying the asymptotic normality of our estimates
reduces to study asymptotic properties of estimators in a
parametric censored nonlinear regression model, such as those
studied by Stute (1999) and Delecroix \& al. (2008). We first recall some elements about the case "$f$ known" (which corresponds to a nonlinear regression setting), and then show that,
under some additional conditions on $\hat{f}$ and on the model,
our estimation of $\theta_0$ is asymptotically equivalent to the one performed in this unreachable parametric model.

\subsection{The case $f$ known}

This case can be studied using the results of Stute (1999) for the WLS approach, or the results of Delecroix \& al. (2008) for the SD approach. We recall some assumptions under which the asymptotic normality of the corresponding estimators is obtained.

\textbf{Assumptions on the model.}
We denote by $\nabla _{\theta }f(x;\theta)$ the vector of partial
derivatives of $(x,\theta)\rightarrow f(\theta'x;\theta)$ with
respect to $\theta$, and $\nabla^2 _{\theta }f$ the corresponding
Hessian matrix.
\begin{assum}
\label{a12} $f(\theta'x;\theta)$ is twice continuously
differentiable with respect to $\theta,$ and $\nabla _{\theta }f$
and $\nabla^2 _{\theta }f$ are bounded as functions of $x$ and
$\theta.$
\end{assum}



\textbf{Assumptions on the censoring. } We need some additional
integrability condition. We first need a moment assumption which is related to the fact that we need to have $E[Y^{*4}]<\infty.$

\begin{assum}
\label{moment} $$\int \frac{y^4dF(y)}{[1-G(y-)]^3}<\infty.$$
\end{assum}
Actually $x$ is not involved in Assumption \ref{moment} as it is assumed to be bounded. Furthermore, in the case $f$ known, this assumption can be weakened, but it will be needed in the case $f$ unknown to obtain uniform consistency rate for $\hat{f}.$ The following assumption is used in Stute
(1995, 1996) to achieve asymptotic normality of Kaplan-Meier
integrals.

\begin{assum}\label{a13}
Let $$C(y)=\int_{-\infty}^{y}\frac{dG(s)}{\{ 1-H(s)\}\{1-G(s)
\}}.$$ Assume that$$\int y C^{1/2}(y)dF_{(X,Y)}(x,y)<\infty.$$
\end{assum}

See Stute (1995) for a full discussion on this kind of assumption.
Using our kernel estimator for estimating the conditional
expectation will lead us to a slightly stronger assumption (see
the appendix section), which is

\begin{assum}
\label{a14}For some $\varepsilon >0$, $$\int y
C^{1/2+\varepsilon}(y)[1-G(y-)]^{-1}dF_{(X,Y)}(x,y)<\infty.$$
\end{assum}

In the following, we will use the (slightly) stronger Assumption
\ref{a14} since it may simplify some proofs (see Lemma \ref{jump2}
and the proof of Theorem \ref{nonlinear}). However, Assumption \ref{a14} could be replaced by Assumption \ref{a13} if we were to use an estimator (not necessarly kernel estimator) which would not require Assumption \ref{a14} to satisfy the proper convergence assumptions. Note that this kind of
assumption is classical in studying regression models with
censored responses. Although it is not mentioned in Burke and Lu
(2005), a similar assumption is implicitly needed to obtain
equation (2.29) of Lai \& al. (1995). In their proof of Lemma A.7
page 199 of Burke \& al. (2005), the authors refer to equation (2.29) page 275 of Lai \& al.
(1995): this only holds under the condition C3 of Lai \& al. (1995) which basically controls the tail behavior of the distributions.

The following Theorem can be deduced from the proof of Theorem 1.2 in Stute (1999) and of Theorem 4 in Delecroix \& al. (2008). However, to make this article self-contained, a short proof of this result is postponed at section \ref{secpost} of the appendix.

\begin{theorem} 
\label{nonlinear}
Define
$$\psi(y,T,\delta)=\left[ \frac{\left( 1-\delta
\right) \mathbf{1}_{T>y}}{1-H\left( T-\right) }-\int \frac{\mathbf{1}%
_{T>y,y>v}dG\left( v\right) }{\left[ 1-H\left( v\right) \right]
^{2}}\right]$$ and let \begin{eqnarray*}
U^{WLS} &=&\frac{\delta \left( T-f\left( \theta _{0}^{\prime
}X;\theta _{0}\right) \right) }{1-G\left( T-\right) }+\int \{
y-f\left( \theta _{0}^{\prime }x;\theta _{0}\right)
\}V(y,T,\delta)
dF_{\left( X,Y\right) }\left( x,y\right) , \\
U^{SD} &=&\left[ \frac{\delta T}{1-G\left( T-\right) }-f\left(
\theta _{0}^{\prime }X;\theta _{0}\right) \right] +\int
yV(y,T,\delta) dF_{\left( X,Y\right) }\left( x,y\right),
\end{eqnarray*}
and let $W^{WLS} =E\left[ \left( U^{WLS}\right) ^{2}\right]$ and $
W^{SD} =E\left[ \left( U^{SD}\right) ^{2}\right].$ Let $M_n$ and $M_{\infty}$ denote respectively either $M_n^{WLS}$ and $M,$ or $M_n^{SD}$ and $M^*.$ We have, under Assumptions \ref{a1} to \ref{a14},
\begin{equation}
M_{n}\left( \theta,f\right) =M_{\infty}\left( \theta ,
f
\right) +O_{P}\left( \frac{\| \theta -\theta _{0}\| }{\sqrt{n}}%
\right) +o_{P}\left( \| \theta -\theta _{0}\| ^{2}\right) +R_{n},
\label{sherman1}
\end{equation}
\begin{equation}
M_{n}\left( \theta,f\right) =\frac{1}{2}%
\left( \theta -\theta _{0}\right) ^{\prime }V\left( \theta -\theta
_{0}\right) +\left( \theta -\theta _{0}\right) ^{\prime }\frac{W_{n}}{\sqrt{n%
}}+o_{P}\left( n^{-1}\right) +R_{n},  \label{sherman2}
\end{equation}
where $R_{n}$ does not depend on $\theta,$ where
\begin{eqnarray*}
V &=&E\left[ \nabla _{\theta }f\left( X;\theta _{0}\right) \nabla
_{\theta }f\left( X;\theta
_{0}\right) ^{\prime }\right],
\end{eqnarray*}
and where
$W_n\Longrightarrow \mathcal{N}(0,W),$ for $W=W^{WLS}$ and $W=W^{SD}$ in the $WLS-$case and $SD-$case respectively.
\end{theorem}

In view of Theorem 1 and 2 of
Sherman (1994), (\ref
{sherman1}) states that, in the case where $f$ is known, $|\hat{\theta}-\theta_0|=O_{P}\left( n^{-1/2}\right) $, while (%
\ref{sherman2}) gives the asymptotic law of $\hat{\theta}$, showing that $n^{1/2}(\hat{\theta}-\theta_0)\Longrightarrow \mathcal{N}(0,V^{-1}WV^{-1}),$ in both WLS and SD cases.

\subsection{The case $f$ unknown}

As $f$ is unknown in the SIM model, we need to add some conditions about the rate of
convergence of $%
\hat{f}$.

\textbf{Assumptions on $f$.} If we evaluate the function $\nabla
_{\theta }f(x;\theta)$ at the point $(x,\theta_0),$ a direct
adaptation of Lemma A.5 of Dominitz and Sherman (2003) shows that
\begin{equation}
\label{gradient} \nabla _{\theta }f(x;\theta_0)=f'(\theta_0'x)\{
x-E[X \mid \theta_0'X=\theta_0'x]\},
\end{equation}
where $f'$ denotes the derivative with respect to $t$ of the
function $f(t;\theta_0).$
\begin{assum}
\label{a16} We assume that the function $f(t;\theta_0)$ is
continuously derivable with respect to $t,$ its derivative is
denoted as $f'$ and is bounded.
\end{assum}

We will also assume some regularity on the model.
\begin{assum}
\label{donsker3} $u\rightarrow f(u;\theta_0)$ where $u$ ranges over $\theta_0'\mathcal{X}$ is assumed to belong to some Donsker class of functions $\mathcal{F}.$
\end{assum}
In our minds, $\mathcal{F}$ will be the class $\mathcal{C}^{1}(\theta_0'\mathcal{X},M),$ that is the class of functions $\phi$ defined on $\theta_0'\mathcal{X}$ and being one time differentiable with $\|\phi\|_{\infty}+\|\phi'\|_{\infty}\leq M$ (see section 2.7 in Van der Vaart and Wellner, 1996). It is important not to impose to much regularity on the regression model, since, as we will see it in Assumption \ref{a15}, $\hat{f}$ will also be required to belong to this class with probability tending to one.


\textbf{Assumptions on }$\hat{f}$.

\begin{assum}
\label{a15}
With probability tending to one, $u\rightarrow \hat{f}(u;\theta_0)\in \mathcal{F}$ where $\mathcal{F}$ is defined in Assumption \ref{donsker3}. Furthermore,
\begin{equation}
\label{c1} \|\nabla_{\theta}\hat{f}-\nabla_{\theta}f\|_{\infty} = o_P(1),
\end{equation}
and, defining $W_i^*=\delta_in^{-1}[1-G(T_i-)]^{-1},$
\begin{align}
\label{c2}  \sup_{\theta\in \Theta_n}\left|\sum_{i=1}^n W_i^*J(\theta_0'X_i)[T_i-f(\theta_0'X_i;\theta_0)][\nabla_{\theta}\hat{f}(X_i;\theta_0)-\nabla_{\theta}f(X_i;\theta_0)]\right| =& o_P(n^{-1/2}), \\
\label{c3}  \sup_{\theta\in \Theta_n}\left|\sum_{i=1}^n W_i^*J(\theta_0'X_i)(\hat{f}(\theta_0'X_i;\theta_0)-f^*(\theta_0'X_i;\theta_0))(\nabla_{\theta}\hat{f}(X_i;\theta)-\nabla_{\theta}f^*(X_i;\theta))\right| =& o_P(n^{-1/2}).
\end{align}
\end{assum}

We can now enounce our asymptotic normality theorem.
\begin{theorem}
\label{an} Under Assumptions \ref{a3} to \ref{a15}, we have
\begin{eqnarray*}
\sqrt{n}\left( \hat{\theta}^{WLS}-\theta _{0}\right) & \Rightarrow
& \mathcal{N}\left( 0,V^{-1}W^{WLS} V^{-1}\right), \\\sqrt{n}\left(
\hat{\theta}^{SD}-\theta _{0}\right) & \Rightarrow &
\mathcal{N}\left( 0,V^{-1}W^{SD} V^{-1}\right).
\end{eqnarray*}
\end{theorem}

\begin{proof}
First apply Proposition \ref{trimming} to obtain that
$J_n(\theta_n'X_i)$ can be replaced by $J(\theta_0'X_i)$ or by
$\mathbf{1}_{f_{\theta}(\theta'X_i)>c/2}$, plus some arbitrary
small terms which will not be mentioned in the following. Moreover, we consider $\theta \in \Theta_n$ which is
an $o_P(1)-$neighborhood of $\theta_0$.

\textbf{Proof for the WLS approach.} Using the representation (\ref{jump}) of the Kaplan-Meier weights,
\begin{eqnarray*}
M_{n}\left( \theta ,\hat{f}\right) &=&M_{n}\left(
\theta ,f\right)-\frac{2}{n}\sum_{i=1}^{n}\frac{\delta _{i}J(\theta_0'X_i)\left(
T_{i}-f\left( \theta
^{\prime }X_{i};\theta \right) \right) }{1-\hat{G}\left( T_{i}-\right) }%
\\ && \times \left[ \hat{f}\left( \theta ^{\prime }X_{i};\theta \right)
-f\left( \theta
^{\prime }X_{i};\theta \right) \right] \\
&&+\frac{1}{n}\sum_{i=1}^{n}\frac{\delta
_{i}J(\theta_0'X_i)}{1-\hat{G}\left( T_{i}-\right) }\left[
\hat{f}\left( \theta ^{\prime }X_{i};\theta \right) -f\left(
\theta
^{\prime }X_{i};\theta \right) \right] ^{2} \\
&=&M_{n}\left( \theta, f\right) -2A_{1n}+B_{1n}.
\end{eqnarray*}
First decompose $A_{1n}$ into four terms,
\begin{eqnarray*}
A_{1n} &=& \frac{1}{n}\sum_{i=1}^{n}\frac{\delta
_{i}J(\theta_0'X_i)\left( T_{i}-f\left( \theta _{0}^{\prime
}X_{i};\theta _{0}\right) \right) }{1-\hat{G}\left( T_{i}-\right)
}\left[ \hat{f}\left( \theta _{0}^{\prime }X_{i};\theta
_{0}\right)
-f\left( \theta _{0}^{\prime }X_{i};\theta _{0}\right) \right] \\
&&+\frac{\delta _{i}J(\theta_0'X_i)\left( f\left( \theta
_{0}^{\prime }X_{i};\theta
_{0}\right) -f\left( \theta ^{\prime }X_{i};\theta \right) \right) }{1-\hat{G%
}\left( T_{i}-\right) }\\
&&\times \left[ \hat{f}\left( \theta ^{\prime }X_{i};\theta
\right) -f\left( \theta ^{\prime }X_{i};\theta \right)
-\hat{f}\left( \theta _{0}^{\prime }X_{i};\theta _{0}\right)
+f\left( \theta _{0}^{\prime }X_{i};\theta
_{0}\right) \right] \\
&&+\frac{\delta _{i}J(\theta_0'X_i)\left( f\left( \theta
_{0}^{\prime }X_{i};\theta
_{0}\right) -f\left( \theta ^{\prime }X_{i};\theta \right) \right) }{1-\hat{G%
}\left( T_{i}-\right) }\left[ \hat{f}\left( \theta _{0}^{\prime
}X_{i};\theta
_{0}\right) -f\left( \theta _{0}^{\prime }X_{i};\theta _{0}\right) \right] \\
&&+\frac{\delta _{i}J(\theta_0'X_i)\left( T_{i}-f\left( \theta
_{0}^{\prime }X_{i};\theta _{0}\right) \right) }{1-\hat{G}\left(
T_{i}-\right)
}\\
&& \times \left[ \hat{f}\left( \theta ^{\prime }X_{i};\theta
\right) -f\left( \theta ^{\prime }X_{i};\theta \right)
-\hat{f}\left( \theta _{0}^{\prime }X_{i};\theta _{0}\right)
+f\left(
\theta _{0}^{\prime }X_{i};\theta _{0}\right) \right] \\
&=&A_{2n}+A_{3n}+A_{4n}+A_{5n}\text{.}
\end{eqnarray*}
$A_{2n}$ does not depend on $\theta $.

For $A_{3n}$, use Assumption \ref{a6} to bound
$f(\theta_0'X;\theta_0)-f(\theta'X;\theta)$ by $M\times
\|\theta-\theta_0\|$ (for some constant $M>0$) using a Taylor expansion. Using a Taylor
expansion, the bracket in $A_{3n}$ can be rewritten as
$(\theta-\theta_0)'[\nabla_{\theta}\hat{f}(X;\tilde{\theta})-\nabla_{\theta}f(X;\tilde{\theta})]$
for some $\tilde{\theta}\in \Theta_n.$ Moreover, using Proposition
\ref{trimming}, we can replace $J(\theta_0'X)$ by
$\mathbf{1}_{\{f_{\tilde{\theta}}(\tilde{\theta}'X)>c/2\}}.$ Hence we have
\begin{eqnarray*}
A_{3n} &\leq & M\|\theta-\theta_0\|^2\sup_{\theta\in \Theta,x\in \mathcal{X}}|\nabla_{\theta}\hat{f}(x;\theta)-\nabla_{\theta}f(x;\theta)|\int d\hat{F}_{(X,Y)}(x,y).
\end{eqnarray*}
The uniform consistency of $\nabla_{\theta}\hat{f}$ in Assumption \ref{a15} shows that $A_{3n}=o_P(\|\theta-\theta_0\|^2).$

For $A_{4n},$ use a second order Taylor expansion and the uniform consistency of $\hat{f}$ to obtain
\begin{eqnarray}
\label{progres}
A_{4n} &=&\frac{1}{n}\sum_{i=1}^{n}\frac{\delta
_{i}J(\theta_0'X_i)\left( \theta -\theta _{0}\right)
}{1-\hat{G}\left( T_{i}-\right) }^{\prime }\nabla _{\theta
}f\left( X_{i};\theta_0\right) \left[ \hat{f}\left( \theta
_{0}^{\prime }X_{i};\theta_0 \right) -f\left( \theta _{0}^{\prime
}X_{i};\theta _{0}\right) \right] \\&&+o_P(\|\theta-\theta_0\|^2),
\end{eqnarray}
In the first term, first replace $G$ by $\hat{G}.$ Using Lemma \ref{jump2} ii) with $\eta=1,$ this introduces a remainder term which is bounded by
$$\frac{O_P(\|\theta-\theta_0\|n^{-1/2})\|\hat{f}-f\|_{\infty}}{n}\sum_{i=1}^{n}\frac{\delta
_{i}J(\theta_0'X_i)C^{1/2+\varepsilon}(T_i-)
}{1-G\left( T_{i}-\right) },$$
where we also used the boundedness of $\nabla_{\theta}f.$ Using the uniform consistency of $\hat{f}$ shows that replacing $\hat{G}$ by $G$ in (\ref{progres}) only arises an $o_P(\|\theta-\theta_0\|n^{-1/2})$ term.
Now, we will use the regularity assumption (\ref{donsker3}) on $f(\cdot;\theta_0).$ If the class $\mathcal{F}$ is Donsker, the class of function $\mathcal{F}'=(\delta,T,X)\rightarrow \delta J(\theta_0'X)[1-G(T_i-)]^{-1}\nabla_{\theta}f(X_i;\theta_0)\mathcal{F}(\theta_0'X_i)$ is Donsker, from a stability property of Donsker classes (see e.g. Van der Vaart and Wellner, 1996). The notation $\mathcal{F}(\theta_0'X_i)$ is used to mention that the functions in $\mathcal{F}$ are evaluated at $\theta_0'X_i.$ Furthermore, for all $\phi\in \mathcal{F}',$
$E[\phi(T_i,\delta_i,X_i)]=0,$ since
$$E\left[\frac{\delta_i\nabla_{\theta}f(X_i;\theta_0)}{1-G(T_i-)}|\theta_0'X_i\right]=E\left[\nabla_{\theta}f(X_i;\theta_0)|\theta_0'X_i\right]=0,$$
from (\ref{gradient}). Hence, using the fact that $\hat{f}(\cdot;\theta_0)\in \mathcal{F}$ with probability tending to one, and the asymptotic equicontinuity property of Donsker classes for $\mathcal{F}'$ (see Van der Vaart and Wellner, 1996), we obtain
$$\frac{1}{n}\sum_{i=1}^{n}\frac{\delta
_{i}J(\theta_0'X_i)\left( \theta -\theta _{0}\right)
}{1-G\left( T_{i}-\right) }^{\prime }\nabla _{\theta
}f\left( X_{i};\theta_0\right) \left[ \hat{f}\left( \theta
_{0}^{\prime }X_{i};\theta_0 \right) -f\left( \theta _{0}^{\prime
}X_{i};\theta _{0}\right) \right]=o_P(\|\theta-\theta_0\|n^{-1/2}),$$
and finally,
$A_{4n}=o_P(\|\theta-\theta_0\|n^{-1/2}).$

Similarly, for $A_{5n}$, a Taylor expansion yields
\begin{eqnarray*}
A_{5n} &=&\frac{\left( \theta -\theta _{0}\right) ^{\prime }}{n}%
\sum_{i=1}^{n}\frac{\delta _{i}J(\theta_0'X_i)\left( T_{i}-f\left(
\theta _{0}^{\prime }X_{i};\theta _{0}\right) \right) }{1-\hat{G}\left(
T_{i}-\right) } \\& & \times \left[ \nabla _{\theta }\hat{f}\left(
X_{i};\tilde{\theta}\right) -\nabla _{\theta }f\left(
X_{i};\tilde{\theta}\right) \right] \\
&=& \frac{\left( \theta -\theta _{0}\right) ^{\prime }}{n}%
\sum_{i=1}^{n}\frac{\delta _{i}J(\theta_0'X_i)\left( T_{i}-f\left(
\theta _{0}^{\prime }X_{i};\theta _{0}\right) \right) }{1-G\left(
T_{i}-\right) }\\
& & \times \left[ \nabla _{\theta }\hat{f}\left(
X_{i};\tilde{\theta}\right) -\nabla _{\theta }f\left(
X_{i};\tilde{\theta}\right) \right]+o_P(\|\theta-\theta_0\|n^{-1/2}),
\end{eqnarray*}
where, as for $A_{4n}$, we replaced $\hat{G}$ by $G$ by using Lemma \ref{jump2} ii) and the uniform consistency of $\nabla_{\theta}\hat{f}.$ Now
We then obtain $A_{5n}=o_P(\|\theta-\theta_0\|n^{-1/2})+o_P(\|\theta-\theta_0\|^2)$ using condition \ref{c3} in Assumption \ref{a15}.

For $B_{1n}$, write
\begin{eqnarray*}
B_{1n} &=&\frac{1}{n}\sum_{i=1}^{n}\frac{\delta
_{i}J(\theta_0'X_i)}{1-\hat{G}\left( T_{i}-\right) }\\
&& \times \left[ \hat{f}\left( \theta ^{\prime }X_{i};\theta
\right) -f\left( \theta ^{\prime }X_{i};\theta \right)
-\hat{f}\left( \theta _{0}^{\prime }X_{i};\theta _{0}\right)
+f\left( \theta _{0}^{\prime }X_{i};\theta
_{0}\right) \right] ^{2} \\
&&+\frac{\delta _{i}J(\theta_0'X_i)}{1-\hat{G}\left( T_{i}-\right)
}\left[ \hat{f}\left( \theta _{0}^{\prime }X_{i};\theta
_{0}\right) -f\left( \theta _{0}^{\prime
}X_{i};\theta _{0}\right) \right] \\
&&+\frac{\delta _{i}J(\theta_0'X_i)}{1-\hat{G}\left( T_{i}-\right)
}\left[ \hat{f}\left( \theta _{0}^{\prime }X_{i};\theta
_{0}\right)
-f\left( \theta _{0}^{\prime }X_{i};\theta _{0}\right) \right] \\
&&\times \left[ \hat{f}\left( \theta ^{\prime }X_{i};\theta
\right) -f\left( \theta ^{\prime }X_{i};\theta \right)
-\hat{f}\left( \theta _{0}^{\prime }X_{i};\theta _{0}\right)
+f\left( \theta _{0}^{\prime }X_{i};\theta _{0}\right) \right]
\end{eqnarray*}
Using a second order Taylor expansion and arguments similar to
those used for $A_{3n},$ we obtain that the first term is of order
$o_{P}\left( \| \theta -\theta _{0}\| ^{2}\right) $. The second
term does not depend on $\theta.$ For the third, a first order
Taylor expansion shows that it is bounded by
$$\|\theta-\theta_0\|\|\nabla_{\theta}\hat{f}-\nabla_{\theta}f\|_{\infty}\sup_{x:J(\theta_0'x)=1}|\hat{f}(\theta_0'x;\theta_0)-f(\theta_0'x;\theta_0)|\int d\hat{F}_{(X,Y)}(x,y).$$ Now condition \ref{c2} in Assumption \ref{a15} shows that this is $o_P(\|\theta-\theta_0\|n^{-1/2}).$

We have just shown that
$$
M_{n}\left( \theta, \hat{f}\right) =M_{n}\left(
\theta , f\right) +o_{P}\left( \frac{\| \theta -\theta
_{0}\| }{\sqrt{n}}\right) +o_{P}\left( \| \theta -\theta _{0}\|
^{2}\right) ,
$$
on a set of probability tending to one. Furthermore, using
(\ref{sherman1}) we deduce $\|
\theta -\theta _{0}\| =O_{P}\left( n^{-1/2}\right) $ from Theorem
1 in Sherman (1994),
and since, from (\ref{sherman2}), on  $O_{P}\left( n^{-1/2}\right)-$neighborhoods of $\theta _{0}$%
,
\[
M_{n}\left( \theta ,f_{\theta }\right) =\frac{1}{2}\left( \theta
-\theta
_{0}\right) ^{\prime }V\left( \theta -\theta _{0}\right) +\frac{1}{\sqrt{n}}%
\left( \theta -\theta _{0}\right) ^{\prime }W^{WLS}_{n}+o_{P}\left(
n^{-1}\right) ,
\]
we can apply Theorem 2 of Sherman to conclude on the asymptotic
law.

\textbf{Proof for }$\phi ^{SD}$. Proceed as for $\phi ^{MC}$, the
only difference is in the fact that $\hat{G}$ does not appear in
the terms where $T$ does not appear at the numerator.
\end{proof}

\section{Simulation study}

\label{secsimul}

In this section, we tried to compare the behavior of our estimator
with the estimator proposed by Burke and Lu (2005) who used the
average derivative technique. We considered three configurations.

\begin{center}
\begin{tabular}{ccccc}
\hline \hline Config 1 & & Config 2 & & Config 3 \\ \hline \hline
$\varepsilon \sim \mathcal{N}(0,2)$ & \vline & $\varepsilon\sim
\mathcal{N}(0,1)$ & \vline & $\varepsilon\sim \mathcal{N}(0,1/16)$ \\
$X\sim \mathcal{U}[-2;2]\otimes \mathcal{U}[-2;2]$ & \vline &
$X\sim \mathcal{U}[0;1]\otimes \mathcal{U}[0;1]$ & \vline & $X\sim
\mathcal{B}(0.6)\otimes \mathcal{U}[-1;1]$ \\
$f(\theta'x;\theta)=1/2(\theta'x)^2+1$ & \vline &
$f(\theta'x;\theta)=\frac{2e^{(0.5\theta'x)}}{0.5+\theta'x}$ & \vline
&
$f(\theta'x;\theta)=1+0.1(\theta'x)^2$ \\
 & \vline & & \vline & $-0.2(\theta'x-1)$ \\
$\theta_0=(1,1)'$ & \vline & $\theta_0=(1,2)'$ & \vline & $\theta_0=(1,2)'$ \\
$C\sim \mathcal{U}[0,\lambda_1]$ & \vline & $C\sim
\mathcal{E}(\lambda_2)$ & \vline & $C\sim \mathcal{E}(\lambda_3)$
\\
\hline \hline
\end{tabular}
\end{center}

The first configuration is used by Burke and Lu (2005) in their
simulation study. Observe that, in this model, (\ref{queue}) does
not hold (this condition (\ref{queue}) is also needed in Burke and
Lu's approach), but it only introduces some asymptotic bias in the
estimation. In the second configuration, there is no such problem
since $C$ is exponential. In the third configuration, we see that
$X$ does not have a Lebesgue density, but $\theta'X$ does. In this
situation, it is expected that the average derivative techniques
does not behave well since it requires that $X$ has a density.

In each configuration, we simulated 1000 samples of different size
$n.$ For each sample, we computed $\hat{\theta}_{WLS},$
$\hat{\theta}_{SD},$ and $\hat{\theta}_{AD}$ which denotes the
average derivative estimate computed from the technique of Burke
and Lu (2005). We then evaluated $\|\hat{\theta}-\theta_0\|^2$ for
each estimate, in order to estimate the Mean Squared Error (MSE)
$E[\|\hat{\theta}-\theta_0\|^2].$ We used different values of the
parameters $\lambda_i$ to modify the proportion of censored
responses ($15\%,$ $30\%,$ and $50\%$ respectively). Results are
presented in the table below.

\begin{center}
\begin{tabular}{ccccccccccc}
\hline \hline Config 1 & & & $n=30$ & & $n=50$ & & $n=100$ \\
\hline \hline $\lambda_1=2.4$ & \vline & $\hat{\theta}^{AD}$ &
$4.8656\times 10^{-2}$ & & $2.6822\times 10^{-2}$ & &
$1.1733\times
10^{-2}$ \\
 & \vline & $\hat{\theta}^{WLS}$ & $1.2814\times 10^{-4}$ & &
 $4.0350\times 10^{-5}$ & & $2.0694\times 10^{-5}$ \\
 & \vline & $\hat{\theta}^{SD}$ & $1.2200\times 10^{-4}$ & & $8.3869\times
 10^{-5}$ &  & $1.3820\times 10^{-5}$ \\
 \hline
 $\lambda_1=1.17$ & \vline & $\hat{\theta}^{AD}$ & $4.5757\times
 10^{-2}$ & & $ 3.3285\times 10^{-2}$ & & $1.8236\times 10^{-2}$ \\
 & \vline & $\hat{\theta}^{WLS}$ &  $1.5713\times 10^{-4}$ & &
 $3.8088\times 10^{-5}$ & & $2.9482 \times 10^{-5}$ \\
 & \vline & $\hat{\theta}^{SD}$ &  $1.6925\times 10^{-4}$ & & $4.0177\times
 10^{-5}$ & &  $1.9924\times 10^{-5}$ \\
 \hline
 $\lambda_1=0.1$ & \vline & $\hat{\theta}^{AD}$ & $1.0102\times
 10^{-1}$& & $7.4870\times 10^{-2}$ & & $5.0438 \times 10^{-2}$
 \\
 & \vline & $\hat{\theta}^{WLS}$ &  $8.3666\times 10^{-4}$ & &
 $1.3010\times 10^{-4}$ & & $3.7669\times 10^{-5}$ \\
 & \vline & $\hat{\theta}^{SD}$ &  $1.2000\times 10^{-3}$ & & $6.7356\times
 10^{-5}$ &  & $2.3650\times 10^{-5}$ \\
 \hline \hline Config 2 & & & $n=30$ & & $n=50$ & & $n=100$ \\
\hline \hline $\lambda_2=0.2$ & \vline & $\hat{\theta}^{AD}$ &
$4.1260\times 10^{-1}$ & & $3.6920\times 10^{-1}$ & &
$3.4151\times 10^{-1}$ \\
 & \vline & $\hat{\theta}^{WLS}$ &  $7.8201\times 10^{-3}$ & &
 $6.5401\times 10^{-3}$ & & $5.8660\times 10^{-3}$ \\
& \vline & $\hat{\theta}^{SD}$ & $1.8296 \times
 10^{-2}$  & & $1.4721\times 10^{-2}$ &  & $1.1034\times 10^{-2}$ \\
 \hline
 $\lambda_2=0.1$ & \vline & $\hat{\theta}^{AD}$ &
$3.5199\times 10^{-1}$ & & $3.3522\times 10^{-1}$ & &
$2.8713\times 10^{-1}$ \\
 & \vline & $\hat{\theta}^{WLS}$ &  $1.2301\times 10^{-2}$ & &
 $7.8301\times 10^{-3}$ & & $7.7180\times 10^{-3}$ \\
 & \vline & $\hat{\theta}^{SD}$ & $2.0822 \times
 10^{-2}$  & & $2.0301\times 10^{-2}$ &  & $1.9741\times 10^{-2}$
 \\\hline
 $\lambda_2=0.05$ & \vline & $\hat{\theta}^{AD}$ &
$1.6238$ & & $1.5553$ & &
$1.5223$ \\
& \vline & $\hat{\theta}^{WLS}$ &  $1.6312\times 10^{-2}$ & &
 $1.5100\times 10^{-2}$ & & $1.2013\times 10^{-2}$ \\
  & \vline & $\hat{\theta}^{SD}$ & $3.0344 \times
 10^{-2}$  & & $2.7057\times 10^{-2}$ &  & $2.2510\times 10^{-2}$ \\
 \hline \hline Config 3 & & & $n=30$ & & $n=50$ & & $n=100$ \\
 \hline \hline $\lambda_3=11$ & \vline & $\hat{\theta}^{AD}$ &
$>10$ & & $>10$ & &
$>10$ \\
& \vline & $\hat{\theta}^{WLS}$ &  $4.1896\times 10^{-4}$ & &
 $3.1530\times 10^{-4}$ & & $1.7453\times 10^{-4}$ \\
  & \vline & $\hat{\theta}^{SD}$ & $4.6218 \times
 10^{-4}$  & & $1.8696\times 10^{-4}$ &  & $1.5286\times 10^{-4}$
 \\\hline
 $\lambda_3=4$ & \vline & $\hat{\theta}^{AD}$ &
$>10$ & & $>10$ & & $>10$ \\
& \vline & $\hat{\theta}^{WLS}$ &  $9.1584\times 10^{-4}$ & &
 $3.3124\times 10^{-4}$ & & $2.8984\times 10^{-4}$ \\
   & \vline & $\hat{\theta}^{SD}$ & $3.4912 \times
 10^{-4}$  & & $2.3344\times 10^{-4}$ &  & $2.2457\times 10^{-4}$
 \\ \hline
 $\lambda_3=2$ & \vline & $\hat{\theta}^{AD}$ &
$>10$ & & $>10$ & & $>10$ \\
& \vline & $\hat{\theta}^{WLS}$ &  $2.0159\times 10^{-2}$ & &
 $1.1431\times 10^{-2}$ & & $2.4111\times 10^{-4}$ \\
    & \vline & $\hat{\theta}^{SD}$ & $9.0591 \times
 10^{-4}$  & & $2.0668\times 10^{-4}$ &  & $1.9921\times 10^{-4}$
 \\ \hline
\end{tabular}
\end{center}

Globally, the performance of the different estimates shrinks when
the proportion of censored responses increases. Performances of
$\hat{\theta}^{WLS}$ and $\hat{\theta}^{SD}$ are globally similar.
In all tested configurations, $\hat{\theta}^{WLS}$ and
$\hat{\theta}^{SD}$ seem to perform better than
$\hat{\theta}^{AD}.$ As expected, in the situation where $X$ does
not have a density, $\hat{\theta}^{AD}$ does not converge.

\section{Appendix}

\label{secappend}

\subsection{Some results on Kaplan-Meier integrals}

\label{secpost}

In this section, we recall some facts on the behavior of
Kaplan-Meier integrals. First part of this section is devoted to
the i.i.d representation of Kaplan-Meier integrals derived by
Stute (1995, 1996), first in the univariate case, then in presence
of covariates. For this, define, for any function $\phi $,
\[
U_{i}\left( \phi \right) = \int
\phi(x,y)\psi(y,T_i,\delta_i)dF_{(X,Y)}(x,y),
\]
where $\psi$ has been defined in Theorem \ref{nonlinear}. It can be
easily shown that $E\left[ U_{i}\left( \phi \right) \right] =0.$
The following Theorem has been derived by Stute (1996).

\begin{theorem}
\label{TCL} Let $\phi $ be a function satisfying
\[
\int \left| \phi \left( x,y\right) \right| C^{1/2}\left( y\right) dF_{\left(
X,Y\right) }\left( x,y\right) <\infty .
\]
Then
\begin{eqnarray*}
\int \phi \left( x,y\right) d\hat{F}_{\left( X,Y\right) }\left( x,y\right)
&=&\frac{1}{n}\sum_{i=1}^{n}\frac{\delta _{i}\phi \left( X_{i},T_{i}\right)
}{1-G\left( T_{i}-\right) } \\
&&+\frac{1}{n}\sum_{i=1}^{n}U_{i}\left( \phi \right) +o_{P}\left(
n^{-1/2}\right) .
\end{eqnarray*}
\end{theorem}

In view of the expression (\ref{jump}) of the jumps of Kaplan-Meier
estimate, this Theorem shows that, asymptotically, these jumps can be
replaced by the ''ideal'' jumps, say $
W_{i}^{*}=n^{-1}\delta _{i}[1-G\left( T_{i}-\right)]^{-1},
$
plus some perturbation that only appears in the study of the variance (since
its expectation is zero). The following lemma gives some additional
precision on the difference between the jumps $W_{in}$ and the ''ideal''
jumps $W_{i}^{*}$.

\begin{lemma}
\label{jump2}
Recall that $\hat{G}$ is the Kaplan-Meier estimator for the distribution of $C,$ $W_{in}=n^{-1}\delta_i[1-\hat{G}(T_i-)]^{-1}$
and $W_i^*=\delta_i[1-G(T_i-)]^{-1},$ and denote by $T_{(n)}$ the largest observation.

\begin{equation}
\sup_{t\leq T_{(n)}}\frac{1-\hat{G}(t-)}{1-G(t-)}=O_{P}(1)\text{
\quad and\quad\ }\sup_{t\leq T_{(n)}}\frac{1-G(t-)}{1-\hat{G}(t-)}%
=O_{P}(1)\,;  \label{zhou_1}
\end{equation}

ii) For all $0\leq \eta \leq 1$
and $\varepsilon
>0$,
\begin{equation}
\left\vert W_{in}-W_i^* \right\vert \leq
W_i^*\{C\left(
T_{i}\right) \}^{\eta[1/2 +\varepsilon] }\times O_{P}\left( n^{-\eta/2}\right) , \label{genial}
\end{equation}%
where the $O_{P}\left( n^{-\eta/2}\right) $ factor does not
depend on $i$.
\end{lemma}

\begin{proof}

\emph{i)} The first part
of (\ref{zhou_1}) follows from Theorem 3.2.4 in Fleming and
Harrington (1991). The second part follows for instance as a
consequence of Theorem 2.2 in Zhou (1991).

\emph{ii)} Fix $\eta >0$ arbitrarily. Since
$\int_{a}^{\tau_H}C^{-1-2\eta}(y)dC(y)<\infty,$ for some $a>0,$
apply Theorem 1 in Gill (1983) to see that%
\begin{equation}
\sup_{y\leq T_{\left( n\right)}}\left[ C\left( y\right) \right]
^{-1/2-\eta }\left\vert Z(y)\right\vert =O_{P}(1), \label{key_ines}
\end{equation}%
where $ Z=\sqrt{n}\{\hat{G}-G\}\{1-G\}^{-1} $ is the Kaplan-Meier
process. Next, the proof can be completed by using the definitions
of $W_{in},$ $W_i^*,$ property (\ref{zhou_1}), and elementary algebra.
\end{proof}

\subsection{Proof of Theorem \ref{nonlinear}}

In this section, we show that the criterion $M_{n}^{WLS}$ and
$M_{n}^{SD}$ satisfy the conditions (\ref{sherman1}) and
(\ref{sherman2}). The same properties can be also shown for the
synthetic data estimators of Leurgans (1987) and Lai \& al.
(1995). More precisions can be found in Delecroix \& al. (2008).
For the sake of simplicity, we only prove it for $M_n^{WLS}$ since the proof for $M_n^{SD}$ uses similar arguments.

\textbf{Proof for }$M_{n}^{WLS}$. Write
\begin{eqnarray}
M_{n}^{WLS}\left( \theta ,f\right) -M\left( \theta \right) &=&
2\int \left( y-f\left( \theta _{0}^{\prime }x;\theta _{0}\right)
\right) \{ f\left( \theta _{0}^{\prime }x;\theta _{0}\right)
-f\left( \theta ^{\prime }x;\theta \right) \} \nonumber \\
&&\times d( \hat{F}_{\left( X,Y\right) }-F_{\left( X,Y\right)
})( x,y)  \nonumber \\
&&+\int \{ f\left( \theta _{0}^{\prime }x;\theta _{0}\right)
-f\left( \theta ^{\prime }x;\theta \right) \} ^{2} \nonumber \\
&&\times d( \hat{F}_{(
X,Y) }-F_{( X,Y) }) ( x,y)  \nonumber \\
&&+\int \left( y-f\left( \theta _{0}^{\prime }x;\theta _{0}\right)
\right) ^{2}d( \hat{F}_{( X,Y) }-F_{( X,Y) }) ( x,y) . \label{MM}
\end{eqnarray}
The last term does not depend on $\theta $. Let $$\chi \left(
x,y\right) =\{ y-f\left( \theta _{0}^{\prime }x;\theta _{0}\right)
\} \nabla _{\theta }f\left( x;\theta _{0}\right) .$$ Using the
derivability Assumption \ref{a12} and Theorem \ref{TCL}, the first
term in the right-hand side of (\ref{MM}) is
\begin{eqnarray}
&&2\left( \theta _{0}-\theta \right) ^{\prime }\int \chi(x,y)
d\left( \hat{F}_{\left( X,Y\right)
}-F_{\left( X,Y\right) }\right) \left( x,y\right) \nonumber \\
&&+2\left( \theta _{0}-\theta \right) ^{\prime } \left[ \int \{
y-f\left( \theta _{0}^{\prime }x;\theta _{0}\right) \} \nabla
_{\theta }^{2}f( x;\tilde{\theta})
\right. \nonumber \\
&& \times \left. d( \hat{F}_{( X,Y) }-F_{( X,Y) }) \left(
x,y\right) \right] \left(
\theta _{0}-\theta \right)  \nonumber \\
&=&2\left( \theta _{0}-\theta \right) ^{\prime } \left\{ \frac{1}{n}%
\sum_{i=1}^{n}\frac{\delta_i
\chi(X_i,T_i)}{1-G(T_i-)}-E\left[\frac{\delta \chi(X,T)}{1-G(T-)}
\right] \right\} \nonumber \\
&& +2\left( \theta _{0}-\theta \right) ^{\prime }\frac{1}{n}%
\sum_{i=1}^{n}U_{i}\left( \chi \right)  +o_{P}\left( \left\|
\theta -\theta _{0}\right\| ^{2}\right), \label{truc}
\end{eqnarray}
where the $o_P$-rate comes from the boundedness of $\nabla
_{\theta }^{2}f$ and consistency of Kaplan-Meier integrals.
Furthermore, the empirical sums in
(\ref{truc}) weakly converge to centered Gaussian variables at rate $O_P(n^{-1/2})$. For the second term in (%
\ref{MM}), rewrite it as
\[
\left( \theta -\theta _{0}\right) ^{\prime } \left[ \int \left[
\nabla _{\theta }f\left( \tilde{\theta}x;\tilde{\theta}\right)
\nabla _{\theta }f\left( \tilde{\theta}x;\tilde{\theta}\right)
^{\prime }\right] d( \hat{F}_{( X,Y) }-F_{( X,Y) }) ( x,y) \right]
\left( \theta -\theta _{0}\right) .
\]
From the boundedness of $\nabla _{\theta }f$ , deduce that this is $%
o_{P}\left( \left\| \theta -\theta _{0}\right\| ^{2}\right) $. We thus
obtained (\ref{sherman1}). To obtain (\ref{sherman2}), use Theorem \ref{TCL}.

\subsection{Properties of $\hat{f}$}

In this section, we derive some properties of $\hat{f}$ defined by
(\ref{ksvcond}), and show that this estimate satisfies Assumptions
\ref{a7} and \ref{a15}. Our approach consists of comparing $\hat{f}$ to the ideal estimator $f^*$ defined as
\begin{equation}
f^{*}\left( \theta ^{\prime }x;\theta \right) =\frac{%
\sum_{i=1}^{n}Y_{i}^{*}K\left( \frac{\theta ^{\prime }X_{i}-\theta
^{\prime }x}{h}\right) }{\sum_{i=1}^{n}K\left( \frac{\theta
^{\prime }X_{i}-\theta ^{\prime }x}{h}\right) }, \label{fetoile}
\end{equation}
that is the estimator based on the true (uncomputable) $Y_{i}^{*}$. Indeed, $f^*$ is a regular kernel estimator based on uncensored variables, and can be studied by traditional nonparametric kernel techniques.

\textbf{Assumptions on the random variables $X'\theta.$}
\begin{assum}
\label{density} For all $\theta \in \Theta,$ $\theta'X$ has a
density which is continuously derivable, with uniformly bounded
derivative.
\end{assum}

\textbf{Assumptions on the kernel function.}
\begin{assum}
\label{kernel}
\begin{itemize}
\item $K$ is symmetric, positive, twice continuously
differentiable function with $K''$ satisfying a Lipschitz
condition.
 \item $\int K(s)ds=1.$
 \item $K$ has a compact support, say $[-1;1]$.
\end{itemize}
\end{assum}
\textbf{Assumptions on the bandwidth.}
\begin{assum}
\label{bandwidth}
\begin{itemize}
\item $nh^8\rightarrow 0.$  \item $nh^{5}[\log(n)]^{-1}=O(1).$
\end{itemize}
\end{assum}

The first Lemma we propose allows us to obtain uniform convergence rates for the ideal estimator $f^*$ as an immediate corollary.

\begin{lemma}
\label{tobecompleted} Let $K$ be a kernel satisfying Assumption
\ref{kernel}. Let $\tilde{K}$ denote either $K$ or its derivative.
Let $Z$ be a random variable with 4th order moment, with
$m(x)=E[Z|X=x]$ twice continuously differentiable, with
derivatives of order 0, 1 and 2 uniformly bounded. Consider, for
$d=0,1,$ and any vectors $x$ and $x'$ in $\mathcal{X},$  
\begin{eqnarray*}
g_n(\theta,x,x',d) &=&
\frac{1}{nh^{1+d}}\sum_{i=1}^n\tilde{K}\left(\frac{\theta'X_i-\theta'x}{h}\right)
\left[\tilde{K}\left(\frac{\theta'X_i-\theta'x'}{h}\right)\right]^dZ_i.
\end{eqnarray*}
We have, for $d=0,1,$
\begin{eqnarray}
\sup_{\theta,x,x'}|g_n(\theta,x,x',d)-E[g_n(\theta,x,x',d)]| &=&
O_P(n^{-1/2}h^{-[d+1]/2}[\log
n]^{1/2}), \label{1}\\
\sup_{\theta,x:f_{\theta}(\theta'x)>c/2}|E[g_n(\theta,x,x,d)]-E[Z|X=x]|
&=& O(h^2), \label{2}\\
\sup_{\theta,x,x'}|E[g_n(\theta,x,x,1)]| &=& O(1). \label{3}
\end{eqnarray}
\end{lemma}

\begin{corollary}
\label{schnell}
Under Assumption \ref{kernel},
\begin{eqnarray*}
\|f^*-f\| &=& O_P(n^{-1/2}h^{-1/2}[\log n]^{1/2}+h^2), \\
\|\nabla_{\theta}f^*-\nabla_{\theta}f\| &=& O_P(n^{-1/2}h^{-3/2}[\log n]^{1/2}+h^2).
\end{eqnarray*}
\end{corollary}

\begin{proof}
For the bias terms (\ref{2}) and (\ref{3}), this can be done by a
classical change of variables, a Taylor expansion, and the fact
that $\int uK(u)du=0$ and $\int u^2K(u)du<\infty.$

For (\ref{1}), first consider
$$g^{M_n}_n(\theta,x,x',d)=\frac{1}{nh^{1+d}}\sum_{i=1}^n\tilde{K}\left(\frac{\theta'X_i-\theta'x}{h}\right)
\left[\tilde{K}\left(\frac{\theta'X_i-\theta'x'}{h}\right)\right]^dZ_i\mathbf{1}_{Z_i\leq
M_n}.$$ We then follow the methodology of Einmahl and Mason
(2005). From Pakes and Pollard (1989), the family of functions
indexed by $(\theta,x,x',h)$ (which has a constant envelope
function),
$$(X,Z)\rightarrow \tilde{K}\left(\frac{\theta'X-\theta'x}{h}\right)
\left[\tilde{K}\left(\frac{\theta'X-\theta'x'}{h}\right)\right]^d\mathbf{1}_{Z\leq
M_n},$$ satisfies the uniform entropy condition of Proposition 1
in Einmahl and Mason (2005) (condition (ii) in their Proposition
1). The other assumptions in their Proposition 1 hold with
$\beta=\sigma=\tilde{C}M,$ for some constant $\tilde{C}$ not
depending on $M.$ We then can apply Talagrand's inequality (see
Einmahl and Mason, 2005, and Talagrand, 1994), with
$\sigma^2_{\mathcal{G}}=n^{-1/2}h^{-[d+1]/2}.$ Take
$M_n=n^{1/2}h^{1/2}.$ It follows from Talagrand's inequality that
$$\sup_{\theta,x,x'}|g^{M_n}_n(\theta,x,x',d)-E[g^{M_n}_n(\theta,x,x',d)]| =
O_P(n^{-1/2}h^{-[d+1]/2}[\log n]^{1/2}).$$ It remains to consider
$g^{M_n}_n-g_n.$ This difference is bounded by
$\tilde{C}n^{-1}h^{-[1+d]}\sum_{i=1}^n |Z_i|\mathbf{1}_{Z_i\geq
M_n}$ for some constant $\tilde{C}.$ This is a sum of positive
quantities, thus we only have to show that its expectation is
$o_P(n^{-1/2}h^{-[d+1]/2}[\log n]^{1/2}).$ For this, apply
H\"{o}lder inequality to bound this expectation by
$h^{-[1+d]}E[Z^4]^{1/4}\mathbb{P}(Z\geq M_n)^{3/4}.$ Moreover,
$\mathbb{P}(Z\geq M_n)\leq E[Z^4]/M_n^4$ from Tschebychev
inequality, and the result follows.
\end{proof}

Proposition \ref{etoile} below ensures that the difference between $\hat{f}$ and $f^*,$ in view of uniform consistency, is asymptotically negligible. Hence Assumption \ref{a7} can be deduced from the uniform consistency of $f^*.$

\begin{proposition}
\label{etoile}Under Assumptions \ref{a14}, \ref{density}, and
Kernel Assumptions \ref{kernel} and \ref{bandwidth}, we have
$\|\hat{f}-f^*\|_{\infty}+\|\nabla_{\theta}\hat{f}-\nabla_{\theta}f^*\|_{\infty}=o_P(1).$
\end{proposition}

\begin{corollary}
Under the Assumptions of Proposition \ref{etoile}, $\hat{f}$ satisfies Assumption \ref{a7} and condition (\ref{c1}) in Assumption \ref{a15}.
\end{corollary}

\begin{proof}
Let $\hat{f}_{\theta'X}(u)=n^{-1}h^{-1}\sum_{i=1}^n K((\theta'X_i-u)/h).$
We have
\begin{eqnarray}\label{fondamental}
\hat{f}(u;\theta)-f^*(u;\theta) &=& \frac{1}{h}\sum_{i=1}^n \frac{[W_{in}-W_i^*]T_iK\left(\frac{\theta'X_i-\theta'x}{h}\right)}{\hat{f}_{\theta'X}(u)}.
\end{eqnarray}
Now, from uniform consistency of kernel density estimator (see, e.g. Einmahl and Mason, 2005),
$$\sup_{x\in \mathcal{X},\theta\in \Theta}|\hat{f}_{\theta'X}(\theta'x)-f_{\theta'X}(\theta'x)|=o_P(1).$$
Using this result on the set $\{f_{\theta'X}(\theta'x)>c>0\}$, and Lemma \ref{jump2} ii) with $\eta$ sufficiently small, we obtain the bound
\begin{equation}
|\hat{f}(\theta'x;\theta)-f(\theta'x;\theta)| \leq  O_P(n^{-\eta/2}h^{-1})\times \sum_{i=1}^n W_i^*T_iC^{\eta(1/2+\varepsilon)}(T_i-)K\left(\frac{\theta'X_i-\theta'x}{h}\right), \label{brune}
\end{equation}
where the $O_P-$rate does not depend on $\theta$ nor $x.$ Recalling the definition of $W_i^*,$ consider the family of functions indexed by $\theta$ and $x,$ $$\{(T,\delta,X)\rightarrow \delta T[1-G(T-)]^{-1}C^{\eta(1/2+\varepsilon)}(T-)K((\theta'X-\theta'x)/h)\}.$$ This family is Euclidean (see Lemma 22 in Nolan and Pollard, 1987) for an enveloppe $\delta T[1-G(T-)]^{-1}C^{\eta(1/2+\varepsilon)}(T-)$ which is, for $\eta=1/2,$ square integrable from Assumption \ref{a14}. Therefore, using the assumptions on the bandwidth,
$$\sup_{x\in \mathcal{X},\theta\in \Theta}|\sum_{i=1}^n W_i^*T_iC^{\eta(1/2+\varepsilon)}(T_i-)K\left(\frac{\theta'X_i-\theta'x}{h}\right)|=O_P(h)+O_P(n^{-1/2}).$$
Finally, back to (\ref{brune}), this shows that $\|\hat{f}-f\|_{\infty}=O_P(n^{-1/4})=o_P(1).$

Similarly, $\|\nabla_{\theta}\hat{f}-\nabla_{\theta}f\|_{\infty}=O_P(n^{-1/4}h^{-1}).$

Now, to prove the corollary, we have to show the uniform consistency of $f^*$ and $\nabla_{\theta}f^*,$ which can be done applying Theorem 2 in Einmahl and Mason (2005).
\end{proof}

The following Proposition allows us to obtain that $\hat{f}$ satisfies conditions (\ref{c2}) and (\ref{c3}) of Assumption \ref{a15}.

\begin{proposition}
\label{burc}
Let $\|\cdot\|_{\Theta_n}$ denote the supremum of the absolute value over $\Theta_n.$ Under the Assumptions of Proposition \ref{etoile}, we have
\begin{align}
h^2\left\|\sum_{i=1}^n W^*_iJ(\theta_0'X_i)(T_i-f(\theta_0'X_i;\theta_0))(\nabla_{\theta}\hat{f}(X_i;\theta)-\nabla_{\theta}f^*(X_i;\theta))\right\|_{\Theta_n} =& O_P(n^{-1}), \label{cc0}\\
h^2\left\|\sum_{i=1}^n W_i^*J(\theta_0'X_i)(\hat{f}(\theta_0'X_i;\theta_0)-f^*(\theta_0'X_i;\theta_0))(\nabla_{\theta}\hat{f}(X_i;\theta)-\nabla_{\theta}f^*(X_i;\theta))\right\|_{\Theta_n} =& O_P(n^{-1}), \label{cc1}\\
\left\|\sum_{i=1}^n W_i^*J(\theta_0'X_i)(\hat{f}(\theta_0'X_i;\theta_0)-f^*(\theta_0'X_i;\theta_0))(\nabla_{\theta}f(X_i;\theta)-\nabla_{\theta}f^*(X_i;\theta))\right\|_{\Theta_n} =& o_P(n^{-1/2}), \label{cc2} \\
\left\|\sum_{i=1}^n W_i^*(f(\theta_0'X_i;\theta_0)-f^*(\theta_0'X_i;\theta_0))(\nabla_{\theta}\hat{f}(X_i;\theta)-\nabla_{\theta}f^*(X_i;\theta))\right\|_{\Theta_n} =& o_P(n^{-1/2}). \label{cc3}
\end{align}
\end{proposition}

\begin{corollary}
\label{final} Under the assumptions of Proposition \ref{etoile}, $\hat{f}$ satisfies conditions (\ref{c2}) and (\ref{c3}) of Assumption \ref{a15}.
\end{corollary}

\begin{proof}[Proof of Corollary \ref{final}]
To prove (\ref{c3}), according to Proposition \ref{burc}, it remains to show that
$$
\sup_{\theta\in \Theta_n}\left|\sum_{i=1}^n W_i^*J(\theta_0'X_i)(f(\theta_0'X_i;\theta_0)-f^*(\theta_0'X_i;\theta_0))(\nabla_{\theta}f(X_i;\theta)-\nabla_{\theta}f^*(X_i;\theta))\right|=o_P(n^{-1/2}),$$
which can be done following the lines of Lemma C2 in Delecroix \& al. (2008). Similarly, Proposition (\ref{burc}) allows to replace $\hat{f}$ by $f^*.$
\end{proof}

\begin{proof}[Proof of Proposition \ref{burc}]
We only prove (\ref{cc1}) and (\ref{cc3}) since the others are similar.

We first prove (\ref{cc1}). This can be done by studying separately the different terms arising by differentiation with respect to $\theta$ in the definition of $\hat{f}.$ We will only study the term coming from the differentiation of the numerator (since the other is similar), that is
\begin{eqnarray*}
\frac{1}{nh^2}\sum_{i,j} \frac{\delta_iJ(\theta_0'X_i)(\hat{f}(\theta_0'X_i;\theta_0)-f^*(\theta_0'X_i;\theta_0)}{1-G(T_i-)} 
K'\left(\frac{\theta'X_i-\theta'X_j}{h}\right)\hat{f}^{-1}_{\theta'X}(\theta'X_i)[W_j^*-W_{jn}]T_j.
\end{eqnarray*}
By bounding $|K'|$ by $\|K'\|_{\infty}$ and using the convergence rate of $f^*,$ it is easily seen that the terms for $i=j$ can be removed from this double sum, arising an $o_P(n^{-1/2})$ term uniform in $\theta.$
Applying (\ref{fondamental}), we then get that the above quantity is, up to an $o_P(n^{-1/2})$term,
\begin{eqnarray*}
\frac{1}{h^3}\sum_{i\neq j,k} W_i^*J(\theta_0'X_i)K\left(\frac{\theta_0'X_i-\theta_0'X_k}{h}\right)[W^*_k-W_{kn}]T_k 
\\
\times K'\left(\frac{\theta'X_i-\theta'X_j}{h}\right)[W_j^*-W_{jn}]T_j\hat{f}_{\theta'X}(\theta'X_i)^{-1}\hat{f}_{\theta_0'X}(\theta_0'X_i)^{-1}.
\end{eqnarray*}
Again, using Lemma \ref{jump2} ii) with $\eta=1/2$, and bounding $K'$ by $\|K'\|_{\infty}$ allows us to remove the terms for $j=k$ and $i=k.$ For the rest of this triple sum, apply Lemma \ref{jump2} ii) with $\eta=1$ and bound $K'$ by $\|K'\|_{\infty}.$ If follows that the left-hand side of (\ref{cc1}) is bounded, uniformly in $\theta,$ by
$$\frac{O_P(n^{-1}h^{-2})}{n}\sum_{i\neq k,j\neq k,i\neq j} W_j^*W_k^*W_i^*C^{1/2+\varepsilon}(T_j-)C^{1/2+\varepsilon}(T_k-)|T_j||T_k|K\left(\frac{\theta_0'X_i-\theta_0'X_k}{h}\right).$$
The last sum as finite expectation (and does not depend on $\theta$) from Assumption \ref{a14}.

For (\ref{cc3}), again, we will consider only the part of $\nabla_{\theta}\hat{f}$ coming from the differentiation of the numerator, this means that we are trying to bound
\begin{equation}\frac{1}{h^2}\sum_{i,j}W_i^*J(\theta_0'X_i)[f(\theta_0'X_i;\theta_0)-f^*(\theta_0'X_i;\theta_0)]K'\left(\frac{\theta'X_i-\theta'X_j}{h}\right)T_j[W_j^*-W_{jn}]\hat{f}_{\theta'X}(\theta'X_i)^{-1}.\label{moche}\end{equation}
First, let $S_{\tau}$ be the double sum deduced from (\ref{moche}) by introducing $\mathbf{1}_{T_j\leq \tau}$ for some $\tau<\tau_H.$ From Gill (1983), $\sup_{t\leq \tau}|\hat{G}(t)-G(t)||1-G(t)|^{-1}=O_P(n^{-1/2}),$ and consequently, $\sup_{j}|W_{j}^*-W_{jn}|\mathbf{1}_{T_j\leq \tau}=O_P(n^{-1/2}).$ Now, using the uniform convergence rate of $f^*$ and bounding $K'$ by $\|K'\|_{\infty}$ shows that $S_{\tau}=o_P(n^{-1/2})$ for any $\tau<\tau_H.$ To obtain a bound for (\ref{moche}), we then have to make $\tau$ tend to $\tau_H.$ For this, we use the same Cramer-Slutsky argument as Stute (1995) in his proof of the Central Limit Theorem under censoring.

Using Lemma \ref{jump2} ii) with $\eta=1$, observe that $$|S_{\tau_H}-S_{\tau}|\leq O_P(n^{-1/2}h^{-2})\|f^*-f\|_{\infty}\frac{1}{h}\sum_{i=1}^nK\left(\frac{\theta_0'X_i-\theta_0'X_j}{h}\right)\mathbf{1}_{T_j\geq \tau}W_j^*C^{1/2+\varepsilon}(T_j-)W_i^*.$$
The last part does not depend on $\theta$ and its expectation tends to zero as $\tau\rightarrow \tau_H,$ while the rest is $O_P(n^{-1/2}),$ using the convergence rate of $f^*$ and the Assumptions \ref{bandwidth}. Then the Cramer-Slutsky argument of Stute allows us to conclude.
\end{proof}

The only condition that still needs to be checked is that $u\rightarrow \hat{f}(u;\theta_0)\in \mathcal{F},$ where $\mathcal{F}$ is defined in Assumption \ref{donsker3}. This can be done if we specify this class of functions. If $\mathcal{F}=\mathcal{C}^1(\theta_0'X,M),$ it suffices to show that $\sup_u|\hat{f}'(u;\theta_0)-f'(u;\theta_0)|=o_P(1),$ which can be done by using the same method as in Proposition \ref{burc} to replace $f$ by $f^*.$

\subsection{Trimming}
In the following proposition, we show that the trimming
$J_n(\theta_n'x)$ can be replaced by $J(\theta_0'x)$ modulo
arbitrary small terms.
\begin{proposition}
\label{trimming} Let, for any function $\phi,$
$$R_n=\frac{1}{n}\sum_{i=1}^{n}\phi(\theta,\hat{G},\hat{f};T_i,\delta_i,X_i)\left[J(\theta_0'X_i)-J_n(\theta_n'X_i)\right].$$
We have $R_n=o_P\left(\frac{1}{n}\sum_{i=1}^{n}\phi(\theta,\hat{G},\hat{f};T_i,\delta_i,X_i)\right)o_P(n^{-1/2})$.
\end{proposition}
\begin{proof}
For any $\delta > 0,$ we have, with probability tending to one,
$$|J(\theta_0'X_i)-J_n(\theta_n'X_i)|\leq \mathbf{1}_{f_{\theta_0'X}(\theta_0'x)\leq c-\delta, \hat{f}_{\theta_n'X}(\theta_n'X) \geq c}+\mathbf{1}_{[\delta;\infty]}(Z_n),$$
where
$Z_n=\sup_{x}|\hat{f}_{\theta_n'X}(\theta_n'X)-f_{\theta_0'X}(\theta_0'x)|\tilde{J}(x).$
As in Delecroix, Hristache, Patilea (2006) page 737-738, we have
$$R_n=o_P\left(\frac{1}{n}\sum_{i=1}^{n}\phi(\theta,\hat{G},\hat{f};T_i,\delta_i,X_i)\right)+\mathbf{1}_{[\delta;\infty]}(Z_n)\times O_P(1).$$
Note that $\mathbb{P}(n^{1/2}Z_n\geq \delta)\leq
\mathbb{P}(Z_n\geq \delta),$ which tends to zero as $\delta$ tends
to zero. \end{proof}

\end{document}